\documentclass[graybox]{svmult}

\usepackage{mathptmx}       % selects Times Roman as basic font
\usepackage{helvet}         % selects Helvetica as sans-serif font
\usepackage{courier}        % selects Courier as typewriter font
\usepackage{type1cm}        % activate if the above 3 fonts are
\usepackage{makeidx}         % allows index generation
\usepackage{multicol}        % used for the two-column index
\usepackage[bottom]{footmisc}% places footnotes at page bottom

\usepackage[english]{babel}
\usepackage{amsmath, amssymb, amsfonts}
\usepackage{array}
\usepackage{graphicx}
\usepackage{float}
\usepackage{theorem}
\usepackage[textheight=0.8\paperheight]{geometry}
\usepackage{blindtext}
\theoremstyle{plain}

\renewcommand{\a}{\alpha}

\renewcommand{\mathcal}{\mathscr}

\newcommand{\PP}{\mathrm{P}}

\newcommand{\eps}{\varepsilon}

\renewcommand{\a}{\alpha}
\renewcommand{\phi}{\varphi}
\renewcommand{\epsilon}{\varepsilon}

\newcommand{\given}{\,|\,}
\newcommand{\Aupper}{\text{A}}
\renewcommand{\E}{\text{E}}

\newcommand{\lb}{\underline}

\renewcommand{\P}{\mathrm{P}}
\renewcommand{\th}{\theta}

\makeindex             % used for the subject index
\begin{document}

\title*{On Bayesian based adaptive confidence sets for linear functionals}
\author{Botond Szab\'o}
\institute{Botond Szab\'o \at Budapest university of Technology, M\H{u}egyetem rkp. 3-9, H-1111 Budapest, Hungary, \email{bszabo@math.bme.hu}}
\maketitle

\abstract{We consider the problem of constructing Bayesian based confidence sets for linear functionals in the inverse Gaussian white noise model. We work with a scale of Gaussian priors indexed by a regularity hyper-parameter and apply the data-driven (slightly modified) marginal likelihood empirical Bayes method for the choice of this hyper-parameter. We show by theory and simulations that the credible sets constructed by this method have sub-optimal behaviour in general. However, by assuming ``self-similarity'' the credible sets have rate-adaptive size and optimal coverage. 
As an application of these results we construct $L_{\infty}$-credible bands for the true functional parameter with adaptive size and optimal coverage under self-similarity constraint.
}

\section{Introduction}

Uncertainty quantification is highly important in statistical inference. Point estimators without confidence statements contain only a limited amount of information. Bayesian techniques provide a natural and computationally advantageous way to quantify uncertainty by producing credible sets, i.e. sets with prescribed (typically $95\%$) posterior probability. In this paper we investigate the validity of such sets from a frequentist perspective. We are interested whether these sets can indeed be used as confidence sets or by doing so one gives a misleading uncertainty quantification, see for instance \cite{Freedman}. We focus in our work on credible sets for linear functionals in nonparametric models and their application to the construction of $L_{\infty}$-credible bands for the functional parameter.

In infinite dimensional models nonparametric priors usually have a tuning- or hyper-parameter controlling the fine details of the prior distribution. In most of the cases the choice of the hyper-parameter in the prior distribution is very influential, incorrect choices can result in sub-optimal behaviour of the posterior. Therefore data driven, adaptive Bayesian techniques are applied in practice to determine the value of the hyper-parameter, overcoming overly strong prior assumptions. The two (perhaps) most well-known Bayesian methods used to achieve adaptive results are the hierarchical Bayes and the empirical Bayes techniques. In our work we focus mainly on the empirical Bayes method, but we conjecture results about the hierarchical Bayes approach as well.

The frequentist properties of adaptive Bayesian credible sets were considered only in a limited number of recent papers; see \cite{SzVZ2, SzVZ3, Belitser, Ray2, Paulo}. The authors of these papers have shown that under a relatively mild and natural assumption on the functional parameter, i.e. the self-similarity condition, the credible sets have good frequentist coverage and in a minimax sense optimal size. However, for non self-similar functions the credible sets provide overconfident, misleading confidence statements. In these papers mainly the $L_2$-norm were considered, which is the natural extension of the finite dimensional Euclidean-norm, but for visualization purposes it is perhaps not the most appropriate choice. In practice usually the posterior credible bands are plotted, which correspond to the $L_{\infty}$-norm. The frequentist properties of $L_{\infty}$-credible bands were investigated in a non-adaptive setting in \cite{CasNickl2, YooGhosal}. Adaptive $L_{\infty}$-credible bands were only considered up to now in the recent work \cite{Ray2}, using techniques developed in \cite{CasNickl2} and \cite{KSzVZ}. In our work we also focus on the construction of $L_{\infty}$-credible bands taking a substantially different and independently developed approach than in \cite{Ray2}. 

In our analysis we consider a sub-class of (possibly non-continuous) linear functionals satisfying the self-similarity condition and construct credible sets over it using empirical Bayes method on a scale of Gaussian priors with varying regularity; see also \cite{SzVZ2} for the application of this family of priors to derive $L_2$ credible sets for the functional parameter. We show that by slightly modifying the empirical Bayes procedure we can construct credible sets with rate adaptive size and good coverage property for the linear functionals of the self-similar functional parameters. However, there exist certain oddly behaving functional parameters (not satisfying the self-similarity assumptions), where the empirical Bayes method provides haphazard and misleading credible sets for the linear functionals of the functional parameter. This result is in itself of independent interest, since until now the frequentist properties of credible sets in semi-parametric problems were mostly investigated in non-adaptive settings, see for instance \cite{bickel2012, Castillo2, Rivoirard, Knapik} and references therein. 

However, perhaps the main contribution of the present paper is the application of the derived results about linear functionals to the analysis of $L_{\infty}$-credible bands. We show that point evaluations of the functional parameter satisfy the self-similarity assumption on linear functionals. Therefore the above described results about credible sets for linear functionals apply also to pointwise credible sets for the functional parameter. Then putting together these pointwise credible sets we arrive at an $L_{\infty}$-credible band, which therefore has again good frequentist properties for self-similar functional parameters.
This technique is essentially different from the one applied in \cite{Ray2} for the construction of $L_{\infty}$-credible bands, where wavelet basis with spike and slab priors were considered and a weak Bernstein-von Mises theorem was proved.

The remainder of the paper is organized as follows. In Section \ref{sec: Main} we introduce the inverse Gaussian white noise model, where we have carried out our analysis. Then in Section \ref{Sec: LinFunc}  we introduce the linear functionals we are interested in (satisfying the self-similarity constraint) and show that the point evaluations of the functional parameter satisfy this property. The construction of the empirical Bayes credible sets are given in Section \ref{Sec: Bayes}. The main results of the paper are formulated in Sections \ref{Sec: Negative} and \ref{Sec: SelfSim}. In Section \ref{sec: sim} we provide a short numerical analysis demonstrating both the positive and negative findings of the paper. The proofs of the main theorems are deferred to Sections \ref{Proof: counter} and \ref{Proof: LinFunc}.

\section{Main result}\label{sec: Main}

Consider the inverse Gaussian white noise model
\begin{align*}
X_t=\int_0^t \mathbb{K}\theta_0(s)ds+\frac{1}{\sqrt{n}}B_t,\qquad t\in[0,1],
\end{align*}
where $B_t$ is the Brownian motion, $1/n$ is the noise level, $X_t$ the observed signal, $\theta_0(\cdot)\in L^2[0,1]$ the unknown function of interest and $\mathbb{K}:L^2[0,1]\mapsto L^2[0,1]$ a given compact, linear, self-adjoint transformation (but we also allow $\mathbb{K}=\mathbb{I}$). From the self-adjoint property of $\mathbb{K}$ follows that its eigenvectors $\phi_i(\cdot):\,[0,1]\mapsto\mathbb{R}$ form an orthogonal basis and the compactness ensures that the corresponding eigenvalues $\kappa_i$ are tending to zero. Hence using series expansion with respect to $\phi_i$ we get the equivalent Gaussian sequence model
\begin{align}
X_i=\kappa_i\theta_{0,i}+\frac{1}{\sqrt{n}}Z_i,\quad\text{for all $i=1,2,...$}\label{def: model}
\end{align}
where $X_i=\langle X_\cdot,\varphi_i(\cdot)\rangle$ and $\theta_{0,i}=\langle \theta_0(\cdot),\varphi_i(\cdot)\rangle$ are the series decomposition coefficients of the observation and the true function, respectively and the random variables $Z_i=\langle B_\cdot,\varphi_i(\cdot)\rangle$ are independent and standard normal distributed. We limit ourselves to the mildly ill-posed inverse problems, where
\begin{align}
C^{-2}i^{-2p}\leq\kappa_i^2\leq C^2 i^{-2p}\label{assump: inverse},
\end{align}
with some fixed non-negative constant $p$ and positive $C$, see \cite{Cavalier2} for the terminology.

Suppose furthermore that the unknown infinite dimensional parameter $\theta_0=(\theta_{0,1},\theta_{0,2},..)$ belongs to a hyper-rectangle
\begin{align}
\Theta^{\beta}(M)=\{\theta:\,\theta_i^2 i^{1+2\beta}\le M\quad\text{for all $i=1,2,...$}\}\label{def: hyperrectangle}
\end{align}
where $\beta$ is the regularity parameter and $M$ is the squared radius of the hyperrectangle. The minimax estimation rate of the full parameter $\theta_0$ is a multiple of $n^{-\beta/(1+2\beta+2p)}$, see \cite{DonohoLiuMcGibbon}.

\subsection{Linear functionals}\label{Sec: LinFunc}
In this paper we focus on the construction of confidence sets for the (possibly unbounded)
linear functionals
\begin{equation}\label{2eq: l}
L\theta = \sum l_i \theta_i,
\end{equation}
where
$l=(l_1,l_2,...)$ is in a self-similar hyper-rectangle $L^{q}_s(R)$, for some $q,R,j_0,K$
\begin{align}
L^{q}_s(R)=\{l\in\ell_2:\, (1/R^2)j^{-1-2q}\leq  \sum_{i=j}^{j+K-1}l_i^2\leq R^2j^{-1-2q}, \text{for all}\, j>j_0\},\label{assump: LinFuncSelf}
\end{align}
where the parameters $j_0$ and $K$ are omitted from the notation. We note that for instance the linear functionals in the form $l_i\asymp i^{-q}$ belong to this hyper-rectangle.

We are particularly interested in the class of non-continuous linear functionals, the point evaluations of the functional parameter $\theta$. For a specific choice of the operator $\mathbb{K}$ all point evaluations on $t\in[0,1]$ belong to $L_s^{-1/2}(R)$; see the next paragraph. Therefore confidence sets for self-similar linear functionals $L\in L_s^{-1/2}(R)$ of the series decomposition coefficients $\theta_0$ also provide us pointwise confidence sets of the function $\theta_0(\cdot)$. Gluing together the (uniform) pointwise confidence sets one arrives to $L_{\infty}$ confidence bands.

In this paragraph we show that point evaluations of the function $\theta_0(\cdot)$ belong to the self-similar class of linear functionals for appropriate choice of the basis. Assume that the eigen-basis of the operator $\mathbb{K}$ is the sine-cosine basis $\phi_i(\cdot)$. The function $\theta_0(\cdot)$ can be given with the help of the trigonometric decomposition
\begin{align}
\theta_0(t)=\sum_{i}\theta_{0,i}\phi_i(t).\label{def: sincosbasis}
\end{align}
Since $\phi_{2i+1}^2(t)+\phi_{2i}^2(t)=\sin^2(i2\pi t)+\cos^2(i2\pi t)=1$ we got that $l_i=\phi_i(t)$ is in $L_s^{-1/2}(2)$ with parameters $j_0=1$ and $K=3$ (since every three consecutive integers contain a pair of $(2i-1,2i)$ for some $i\in\mathbb{N}$). 

\subsection{Bayesian approach}\label{Sec: Bayes}

In the Bayesian framework to make inference about the unknown sequence $\theta_0$ we endow it with a prior distribution. In our analysis we work with the infinite dimensional Gaussian distribution
\begin{align}
\Pi_{\alpha}=\bigotimes _{i=1}^{\infty}N(0,i^{-1-2\alpha}),\label{prior}
\end{align}
where the parameter $\alpha>0$ denotes the regularity level of the prior distribution. One can easily compute the corresponding posterior distribution
\begin{align}
\Pi_\alpha(\cdot\given X) = \bigotimes_{i=1}^{\infty}{N}\Big(\frac{n\kappa_i^{-1}}{i^{1+2\alpha}\kappa_i^{-2} + n}X_i,
\frac{\kappa_i^{-2}}{i^{1+2\alpha}\kappa_i^{-2}+n}\Big).\label{PostDist}
\end{align}
Furthermore by combining and slightly extending the results of \cite{Knapik} and \cite{Castillo} one can see that the choice $\alpha=\beta$ leads to posterior contraction rate $n^{-\frac{\beta}{1+2\beta+2p}}$ for $\theta_0\in\Theta^{\beta}(M)$, while other choices of the parameter $\alpha$ provide sub-optimal contraction rates.

In this paper, however, we are interested in the posterior distribution of the linear functional $L\theta$.
From Proposition 3.2 of \cite{Knapik} follows that the posterior distribution of the linear functionals $L\theta$ (assuming measurability with respect to the prior $\Pi_{\a}$) takes the form
\begin{align}
\Pi_{\a}^{L}(\cdot|X)=N\Big(\sum_i \frac{nl_i \kappa_i^{-1}}{i^{-1-2\a}\kappa_i^{-2}+n}X_i,\sum_i \frac{l_i^2\kappa_i^{-2}}{i^{1+2\a}\kappa_i^{-2}+n}\Big).\label{eq: postlin}
\end{align}
Furthermore it was also shown in Section 5 of \cite{Knapik} that the optimal choice of the hyper-parameter $\a$ is not $\beta$, but rather $\beta-1/2$. The resulting optimal rate is of the order $n^{-(\beta+q)/(2\beta+2p)}\vee n^{-1/2}$; see \cite{Donoho, Donoho2}. Note that in case $q\geq p$ the smoothness of the linear functional compensates for the degree of ill-posedness and we get a regular problem with contraction rate $n^{-1/2}$. However, in our work we focus on the (from the construction of credible bands point of view) more interesting case $q<p$.

Since the regularity parameter $\beta$ of the infinite sequence $\theta_0$ is usually not available one has to use data-driven method to choose $\alpha$, which we will refer from now on as the hyper-parameter of the prior. Following \cite{KSzVZ} and \cite{SzVZ2} we select a value for $\alpha$ with the marginal likelihood empirical Bayes method, i.e. we select the maximizer of
\begin{align}
\hat\alpha_n=\arg\max_{\a\in[0,\Aupper]}\ell_n(\a),\label{def: hatA}
\end{align}
where $A$ is some arbitrary large, fixed constant, and $\ell_n$ denotes the corresponding log-likelihood for $\alpha$ (relative to an infinite
product of $N(0,1/n)$-distributions)
\begin{align}
\ell_n(\alpha)=-\frac{1}{2}\sum_{i=1}^{\infty}\Big( \log\Big(1+\frac{n}{i^{1+2\alpha}\kappa_i^{-2}}\Big)-\frac{n^2} {i^{1+2\alpha}\kappa_i^{-2}+n}X_i^2 \Big).\label{fn(alpha)}
\end{align}
Then the {\em empirical Bayes posterior} for the functional parameter $\theta$ is defined as  $\Pi_{\hat\alpha_n}(\cdot | X)$ obtained by
substituting $\hat\alpha_n$ for $\alpha$ in the posterior distribution (\ref{PostDist}), i.e.\
\[
\Pi_{\hat\alpha_n}(B | X) = \Pi_{\alpha}(B | X) \Big|_{\alpha = \hat\alpha_n}
\]
for measurable subsets $B \subset \ell^2$. Slightly adapting the proof of Theorem 2.3 in \cite{KSzVZ} we can get that the posterior distribution of the functional parameter $\theta$ achieves the corresponding minimax contraction rate up to a logarithmic factor.

As conjectured in page 2367 of \cite{Knapik} and Section 2.3 of \cite{KSzVZ} this suggests that the present procedure is sub-optimal for the linear functional $L\theta_0$, since adaptation for the full parameter $\theta_0$ and its linear functionals $L\theta_0$ is not possible simultaneously in this setting. In view of the findings in the non-adaptive case \cite{Knapik}
we might expect, however, that we can slightly alter the procedures to deal with linear functionals.
For instance, it is natural to expect that the empirical Bayes posterior for linear functionals $L\theta$, given in $\eqref{2eq: l}$,
\begin{align}
\Pi^{L}_{\hat\a_n-1/2}(\cdot\given X)=
\Pi^{L}_{\a}(\cdot\given X)\Big|_{\a=\hat\a_n-1/2}\label{eq: ebpostlin}
\end{align}
yields optimal rates. In the present paper we work with this data-driven choice of the hyper-parameter and investigate the frequentist properties of Bayesian credible sets constructed from the posterior $\eqref{eq: ebpostlin}$.

For fixed hyper-parameter $\a$ the posterior $\eqref{eq: postlin}$ is a one dimensional Gaussian distribution hence a natural choice of the credible set is the interval
\begin{align}
\hat{C}_{n,\a}=[ \widehat{L\theta}_{n,\a}- \zeta_{\gamma/2}s_n(\a),  \widehat{L\theta}_{n,\a}+ \zeta_{\gamma/2}s_n(\a)],\label{def: cred0}
\end{align}
where $ \widehat{L\theta}_{n,\a}$ is the posterior mean, $s_n^2(\a)$ the posterior variance given in $\eqref{eq: postlin}$, and $\zeta_{\gamma}$ is the $1-\gamma$-quantile of the standard normal distribution. We note that the mean $ \widehat{L\theta}_{n,\a}$ of the posterior distribution of the linear functional  is exactly the linear functional $L$ of the posterior mean of the full parameter $\hat\theta_{n,\a}$ given in $\eqref{PostDist}$. One can easily see that the preceding interval accumulates $1-\gamma$ fraction of the posterior mass. Then the empirical Bayes credible sets are obtained by replacing $\a$ with the data-driven choice $\hat\a_n-1/2$ in $\eqref{def: cred0}$. We introduce some additional flexibility, by allowing the blow up of the preceding interval with a constant factor $D>0$, i.e.
\begin{align}
\hat{C}_n^{L}(D)=[ \widehat{L\theta}_{n,\hat\a_n-1/2}- D\zeta_{\gamma/2}s_n(\hat\a_n-1/2),  \widehat{L\theta}_{n,\hat\a_n-1/2}+ D\zeta_{\gamma/2}s_n(\hat\a_n-1/2)].\label{eq: ebcrediblelin}
\end{align}

\vspace{-0.5cm}
\subsection{Negative results}\label{Sec: Negative}

However, matters seem to be more delicate than one would expect from the direct case. For certain oddly behaving functions the posterior distribution $\eqref{eq: ebpostlin}$ achieves only sub-optimal contraction rates. Furthermore the credible sets $\eqref{eq: ebcrediblelin}$  have also coverage tending to zero.

\begin{theorem}\label{2Lemma: CounterExample}
Let $n_j$ be positive integers such that $n_1\geq2$ and $n_j\geq n_{j-1}^4$ for every $j$.
and let $K > 0$.  Let
$\theta_0=(\theta_{0,1},\theta_{0,2},...)$ be such that
\begin{align}
\theta_{0,i}=\begin{cases}
K n_j^{-\frac{1/2+\beta}{1+2\beta+2p}},  & \text{if $n_j^{\frac{1}{1+2\beta+2p}}\leq i<2n_j^{\frac{1}{1+2\beta+2p}}$}, \qquad j=1,2,\ldots, \\
0,  & \text{else},\label{2eq: Counter}
\end{cases}
\end{align}
for some positive constants $\beta,p$. Then the constant $K > 0$ can be chosen such that the coverage of the credible set tends to zero for every $q\in\mathbb{R}$, $D,R>0$ and $L\in L_s^{q}(R)$
\begin{align*}
\P_{\theta_0}(\theta_0\in\hat{C}_{n_j}(D))\rightarrow 0.
\end{align*}
Furthermore the posterior distribution attains sub-optimal contraction rate
\begin{align}
\Pi_{\hat\a_{n_j}-1/2}^{L}(L\theta:\, |L\theta_0-L\theta|\geq m n_j^{-(\beta+q)/(1+2\beta+2p)}|X) \stackrel{\P_{\theta_0}}{\rightarrow} 1,
\end{align}
as $j \to \infty$ for a positive, small enough constant $m$ and linear functional $L$ satisfying $\eqref{assump: LinFuncSelf}$.
\end{theorem}

The proof of the Theorem is given in Section \ref{Proof: counter}. The sub-optimal contraction rate of the posterior distribution and the bad coverage property of the credible sets are due to the mismatch of the underlying loss functions. In the empirical Bayes method the hyper-parameter $\a$ is chosen to maximize the marginal likelihood function. This method is related to minimizing the Kullback-Leibler divergence between the marginal Bayesian likelihood function and the true likelihood function. At the same time the evaluation of the posterior distribution is given with respect to some linear functional $L$ of the functional parameter $\theta_0$. Optimal contraction rate and good coverage follows from optimal bias-variance trade-off. However, the likelihood based empirical Bayes method intends to minimize the Kullback-Leibler divergence, which is not an appropriate approach in general for balancing out the bias and variance terms. Therefore the empirical Bayes (and we believe that also the hierarchical Bayes method) leads to sub-optimal rate and poor coverage.

\vspace{-0.5cm}
\subsection{Self-similarity}\label{Sec: SelfSim}

To solve this problem we can introduce some additional constraint on the regularity class $\Theta^{\beta}(M)$. The notion of self-similarity originates from the frequentist literature \cite{PicTri,GineNickl,Bull,CCK14, NicklSz}, and was adapted in the Bayesian literature \cite{SzVZ2,Ray2, Belitser}. We call a series in the hyper-rectangle $\theta_0\in \Theta^{\beta}(M)$ \emph{self-similar} if it satisfies
\begin{align}
\sum_{i=N}^{\rho N}\theta_{0,i}^2\geq \eps M N^{-2\beta}, \quad\forall  N\geq N_0, \label{prop: SelfSim}
\end{align}
where $N_0, \eps$ and $\rho$ are some fixed positive constants. Furthermore we denote the class of functions satisfying the  self-similar constraint by $\Theta_s^\beta(M)$, where we omit the parameters $N_0,\eps,\rho$ from the notation. We denote by $\Theta_s(M)$ the collection of self-similar functions with regularity in a compact interval of regularity parameters $\beta\in[\beta_{\min},\beta_{\max}]$
\begin{align}
\Theta_s(M)=\cup_{\beta\in[\beta_{\min},\beta_{\max}]}\Theta_s^{\beta}(M).
\end{align}
Here we omit again in the notation the dependence on $\beta_{\min}$ and $\beta_{\max}$ and assume that $\beta_{\min}>-q$, else $L\theta_0$ would be infinite.

We show that uniformly over $L\in L^{q}_s(R)$ and $\theta_0\in\Theta_s(M)$ the coverage of credible sets $\hat{C}^{L}_n(D)$ for the linear functionals $L\theta_0$ tends to one. Furthermore we prove that the size of the credible sets achieves the corresponding minimax contraction rate.

\begin{theorem}\label{thm: CovLin}
There exists a  large enough positive constant $D$ such that the empirical Bayes credible sets $\hat{C}_n^L(D)$ have honest asymptotic coverage one over the self-similar linear functionals $L\in L_s^{q}(R)$ of the functional parameter $\theta_0$ satisfying $\eqref{assump: LinFuncSelf}$, i.e.
\begin{align}
\inf_{\theta_0\in\Theta_s(M)}P_{\theta_0}\big(L\theta_0\in \hat{C}_n^{L}(D),\,\forall L\in L^{q}_s(R) \big)\rightarrow1.\label{eq: LinPostCov}
\end{align}
Furthermore the radius of the credible sets is rate adaptive, i.e. there exists a positive constant $C_1>0$ such that for all $\beta\in(q,\beta_{\max}]$ we have 
\begin{align}
\inf_{\theta_0\in \Theta^{\beta}_s(M)}\P_{\theta_0}\big(s_n(\hat{\alpha}_n-1/2)\leq C_1 n^{-\frac{\beta+q}{2\beta+2p}},\forall L\in L_s^{q}(R)\big)\rightarrow1.\label{eq: Linadaptivity}
\end{align}
\end{theorem}

We defer the proof to Section \ref{Proof: LinFunc}. The credible band on $[0,1]$ can be constructed with the help of the linear functionals $\phi_i(t)$ introduced in $\eqref{def: sincosbasis}$, i.e. the point evaluations of the basis $\phi_i(\cdot)$ at $t\in[0,1]$. Following from its definition $\eqref{eq: ebcrediblelin}$ the credible band takes the form
\begin{align}
[\hat\theta_n(t)-D\zeta_{\gamma/2}s_n(t,\hat\a_n-1/2), \hat\theta_n(t)+D\zeta_{\gamma/2}s_n(t,\hat\a_n-1/2)], \quad t\in [0,1],\label{eq: CredBand}
\end{align}
where $\hat\theta_n(t)$ is the posterior mean and $s_n^2(t,\hat\a_n-1/2)$ is the posterior variance for $\a=\hat\a_n-1/2$ given in $\eqref{eq: ebpostlin}$ belonging to the linear functional $L=(l_i)_{i\geq1}=(\phi_i(t))_{i\geq1}$. By combining Theorem \ref{thm: CovLin} and the argument given in the last paragraph of Section \ref{Sec: LinFunc} we get that the credible band 
$\eqref{eq: CredBand}$ has honest coverage and rate adaptive size.

\begin{corollary}\label{cor: CovLin}
Assume that the eigen-vectors $\phi_i(\cdot)$ of the linear operator $\mathbb{K}$ form the sine-cosine basis. Then there exists a constant $D$ such that the empirical Bayes credible bands, given in $\eqref{eq: CredBand}$, have honest asymptotic coverage one
\begin{align}
\inf_{\theta_0\in\Theta_s(M)} P_{\theta_0}(|\theta_0(t)-\hat\theta_{n,\hat\a_n-1/2}(t)|\leq D\zeta_{\gamma/2}s_n(t,\hat\a_n-1/2),\,\forall t\in[0,1] \big))\rightarrow1.\label{eq: LinPostCov2}
\end{align}
Furthermore the size of the credible band is rate optimal in a minimax sense, i.e. there exists a $C_1>0$ such that for all $\beta\in(1/2,\beta_{\max}]$
\begin{align}
\inf_{\theta_0\in\Theta_s^{\beta}(M)}\P_{\theta_0}\big(s_{n}(t,\hat{\alpha}_n-1/2)\leq C_1 n^{-\frac{\beta-1/2}{2\beta+2p}}\big)\rightarrow1.\label{eq: Linadaptivity2}
\end{align}
\end{corollary}

\vspace{-0.5cm}
\section{Simulation study}
\label{sec: sim}
We investigate our new empirical Bayes method in an example. We consider the model $\eqref{def: model}$ with $\mathbb{K}=\mathbb{I}$ (the identity operator) and work with the sine-cosine basis on $[0,1]$, i.e. $\phi_1(t)=1$, $\phi_{2i}(t)=\sqrt{2}\cos(2\pi it)$, $\phi_{2i+1}(t)=\sqrt{2}\sin(2\pi it)$ for $t\in[0,1]$. 

\begin{figure}[htbp]
  \centering
 \includegraphics[scale=1.05]{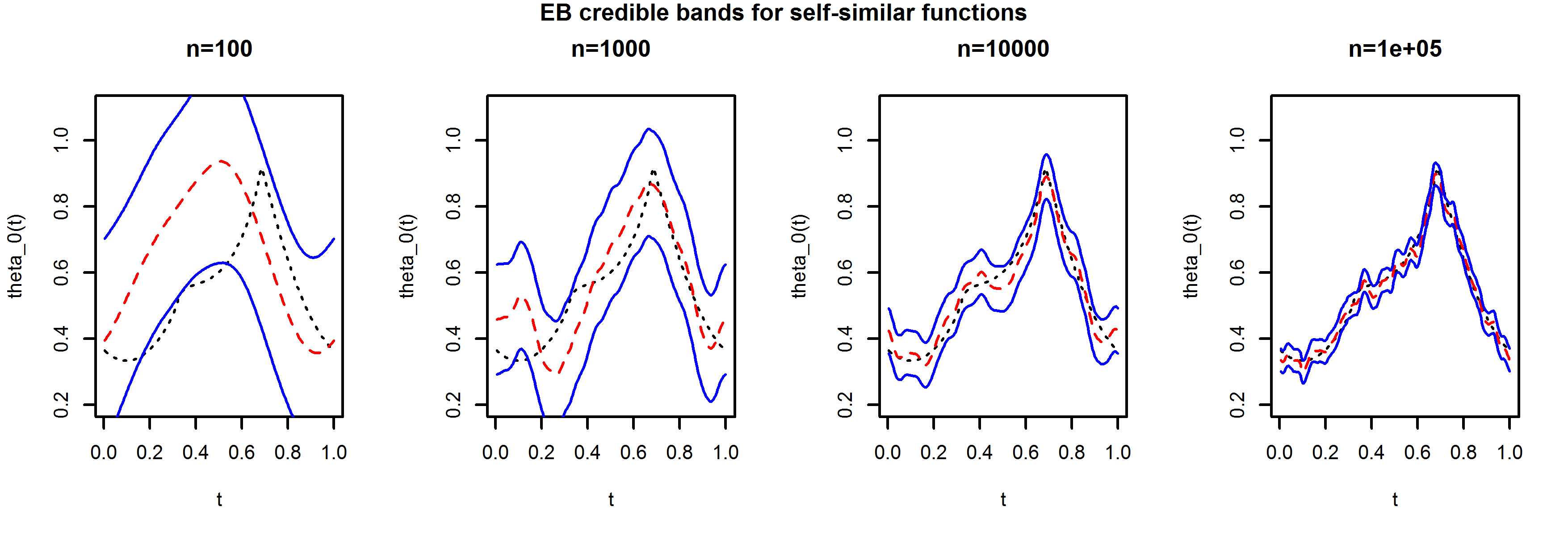}
 \caption{Empirical Bayes credible bands for a self-similar function. The true function is drawn in pointed black, the posterior mean in dashed red and the credible bands in blue. From left to right we have $n=100, 10^3, 10^4$ and $10^5$.}
  \label{figure: SelfSim}
\end{figure}

First we illustrate that for self-similar functions our method provides reasonable and trustworthy credible sets which could be used as confidence bands. We define the true function $\theta_0(t)$ with the help of its sine-cosine basis coefficients, i.e. we take $\theta_{0,i}=i^{-2}\cos(i)$:
\begin{align*}
\theta_0(t)=\cos(1)+\sqrt{2}\sum_{i=1}^{\infty}(2i)^{-2}\cos(2i) \cos(2\pi i t)+\sqrt{2}\sum_{i=1}^{\infty}(2i+1)^{-2}\cos(2i+1) \sin(2\pi i t).
\end{align*}
For computational convenience we work only with the first $10^3$ Fourier coefficients of the true function.
 We simulate data from the corresponding distribution with noise level $n=100, 10^3, 10^4$ and $10^5$. Figure \ref{figure: SelfSim} shows the true function in pointed black, the posterior mean in dashed red and the $95\%$ credible bands (without blowing it up by a constant factor $D$) in blue. One can see that for every noise level $n$ the credible band has good coverage and is concentrating around the truth as $n$ increases, confirming the results of Corollary \ref{cor: CovLin}.

\begin{figure}[htbp]
  \centering
 \includegraphics[scale=1.05]{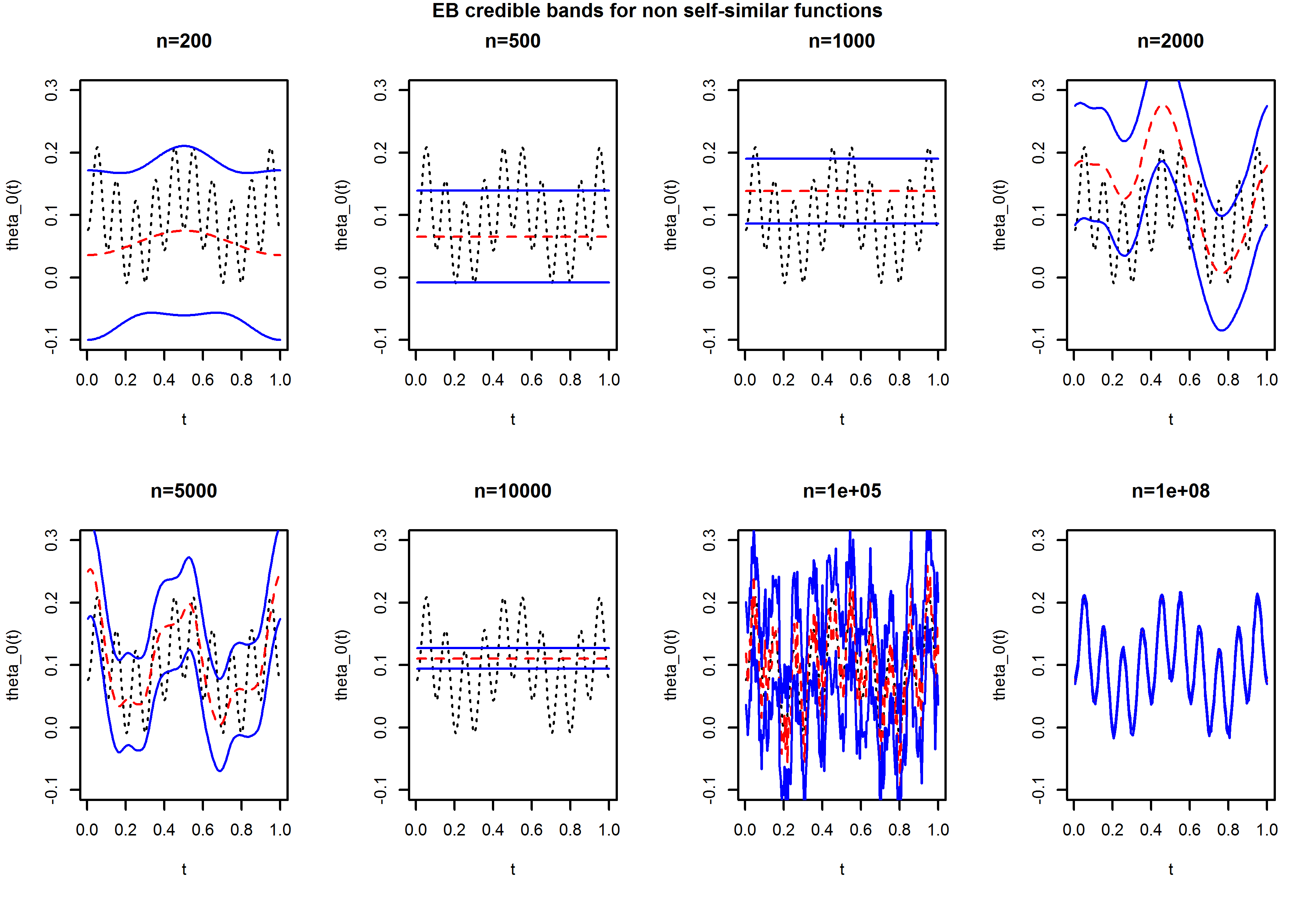}
 \caption{Empirical Bayes credible bands for a non-self-similar function. The true function is drawn in pointed black, the posterior mean in dashed red and the credible bands in blue. From left to right and top to bottom we have  $n=200,500,10^3,2*10^3,5*10^3,10^4,10^5$ and $10^{8}$.}
  \label{figure: CounterExample}
\end{figure}

%\vspace{-1cm}
To illustrate the negative result derived in Theorem \ref{2Lemma: CounterExample} (for the point evaluation linear functionals) we consider a non-self-similar function $\theta_0(t)$ defined by its series decomposition coefficients with respect to the sine-cosine basis. We take the coefficients to be $\theta_{0,1}=1/10,\,\theta_{0,4}=1/30,\,\theta_{0,20}=-1/20,\, \theta_{0,i}=i^{-3/2}$ if $2^{4^j}<i\leq 2*2^{4^j}$ for  $j\geq2$, and $0$ otherwise:
\begin{align*}
\theta_0(t)=0.1+\frac{\sqrt{2}}{30}\cos(4\pi t)-\frac{\sqrt{2}}{20}\cos(20\pi t)+\sum_{j=2}^{\infty} \Big(\sqrt{2}\sum_{i=2^{4^j-1}+1}^{2^{4^j}}(2i)^{-3/2}\cos(2\pi it)+(2i+1)^{-3/2}\sin(2\pi it)\Big).
\end{align*}
For simplicity we consider again only the first $10^3$ Fourier coefficient of the true function. Then we simulate data from the corresponding distribution with various noise levels $n=200,500,10^3,2*10^3,5*10^3,10^4,10^5$ and $10^{8}$. In Figure \ref{figure: CounterExample} we plotted the $95\%$ $L_{\infty}$-credible bands with blue lines, the posterior mean with dashed red line and the true function with pointed black line. One can see that for multiple noise levels we have overly confident, too narrow credible bands 
($n=500,10^3,2*10^3,10^4$), while for other values of the noise levels $n$ we have good coverage ($n=200,5*10^3,5*10^4,10^8$). These periodicity between the good and the bad coverage of the credible sets continues as $n$ increases (but to see it we have to zoom in into the picture).

\vspace{-1cm}

\section{Proof of Theorem \ref{2Lemma: CounterExample}}\label{Proof: counter}
Following \cite{SzVZ}, \cite{KSzVZ} and \cite{SzVZ2} we introduce the notation
\begin{align*}
h_n(\a;\th_0)= \frac{1+2\a+2p}{n^{1/(1+2\a+2p)}\log n}
\sum_{i=1}^{\infty}\frac{n^2i^{1+2\a}(\log i)\th_{0,i}^2}
{(i^{1+2\a+2p}+n)^2}, \qquad \a\ge 0,
\end{align*}
and define
\begin{align*}
\underline{\a}_n(\th_0)
&=\inf\{\a\in [0,\Aupper]: h_n(\a;\th_0)\geq 1/(16C^8)\},\\
\overline{\a}_n(\th_0)
&=\sup\{\a\in [0,\Aupper]: h_n(\a;\th_0)\leq 8C^8 \},
\end{align*}
where the parameter $\Aupper$ was introduced in $\eqref{def: hatA}$.
From the proof of Theorem 5.1 of \cite{SzVZ2} one can see that
\begin{align}
\inf_{\theta_0\in \ell_2}\P(\lb\a_n\leq\hat\a_n)\rightarrow 1.\label{eq: LBhatA}
\end{align}

Furthermore, let us introduce the notations
\begin{align*}
B_n^{L}(\a)=|\E_{\theta_0}\widehat{L\theta}_{\a}-L\theta_{0}|\quad\text{and}\quad V_n(\a)=|\widehat{L\theta}_\a-\E_{\theta_0}\widehat{L\theta}_{\a}|,
\end{align*}
where $\widehat{L\theta}_{\a}$ denotes the posterior mean of the linear functional $L\theta$ for a fixed hyper-parameter $\a>0$.
Similarly to the proof of Theorem 3.1 of \cite{SzVZ2} we have following from the triangle inequality that $\theta_0\in\hat{C}_n^{L}(D)$ implies
$B_n^{L}(\hat\a_n-1/2)\leq V_n(\hat\a_n-1/2)+D\zeta_{\gamma/2}s_n(\hat\a_n-1/2)$. Therefore following from the convergence $\eqref{eq: LBhatA}$ we have
\begin{align}
\P_{\th_0}(L\th_0\in \hat{C}^{L}_n(D))
\leq \P_{\theta_0}\Big(\inf_{\alpha\geq\underline\alpha_n-1/2}B_n(\a)\leq D\sup_{\alpha\geq\underline\alpha_n-1/2}\big[D\zeta_{\frac{\gamma}{2}}s_{n}(\alpha)+V_n(\a)\big]\Big)+o(1).
\label{eq: coverCounter}
\end{align}

From the proof of Theorem 3.1 of \cite{SzVZ2} follows that $\underline\a_{n_j}>\beta+1/2$ for $j$ large enough. Following from the proof of Theorem 5.3 of \cite{Knapik} we get that both $s_{n_j}(\a)$ and $V_{n_j}(\a)$ are bounded from above by a multiple of $n_j^{-(1/2+\beta+q)/(1+2\beta+2p)}\ll n_j^{-(\beta+q)/(2\beta+2p)}$with probability tending to one.

Furthermore, for fixed hyper-parameter $\a$ the bias corresponding to the posterior mean $\eqref{eq: ebpostlin}$ is
$$B_n^{L}(\a)=|\sum_i\frac{l_i \theta_{0,i}}{1+n_ji^{-1-2\a}\kappa_{i}^{2}}|.$$
Note that following from the definition of $L_s^{q}(R)$ given in $\eqref{assump: LinFuncSelf}$ we have that 
$$\sum_{i=j}^{j+K-1}|l_i|\geq \max_{i\in \{j,j+1,...,j+K-1\}}|l_i|\geq j^{-1/2-q}/(RK). $$
Then for $\a\geq\underline\a_{n_j}-1/2\geq\beta$ the squared bias $B_n^2(\a)$ corresponding to the sequence $\eqref{2eq: Counter}$  can be bounded from below by
\begin{align}
 \Big(\sum_{i=n_j^{\frac{1}{1+2\beta+2p}}}^{2n_j^{\frac{1}{1+2\beta+2p}}}\frac{C_0^{-1}|l_i|K n_j^{-\frac{1/2+\beta}{1+2\beta+2p}}}{1+n_ji^{-1-2\a}\kappa_i^{2}}\Big)^2
&\gtrsim n_j^{-\frac{1/2+\beta}{1+2\beta+2p}}
\sum_{i=n_j^{1/(1+2\beta+2p)}/K}^{2n_j^{1/(1+2\beta+2p)}/K-1}
\sum_{j=iK+1}^{(i+1)K}|l_i|\nonumber\\
&\gtrsim n_j^{-\frac{1/2+\beta}{1+2\beta+2p}}\sum_{i=n_j^{1/(1+2\beta+2p)}/K}^{2n_j^{1/(1+2\beta+2p)}/K-1}i^{-1/2-q},\label{eq: LBforBias}
\end{align}
which is further bounded from below by $ n_j^{-(2\beta+2q)/(1+2\alpha+2p)}\gg n_j^{-2(\beta+q)/(2\beta+2p)}$. Therefore the probability on the right hand side of the inequality $\eqref{eq: coverCounter}$ tends to zero.

Finally we note that following from the sub-optimal order of the bias term $\eqref{eq: LBforBias}$ the posterior distribution achieves sub-optimal contraction rate around the true value $L\theta_0$.

\vspace{-0.5cm}
\section{Proof of Theorem \ref{thm: CovLin}}\label{Proof: LinFunc}
First we note that following from the inequalities $(6.9)$ and $(6.10)$ of \cite{SzVZ2} we have for all $\beta\in[\beta_{\min},\beta_{\max}]$  that
\begin{align}
\beta-K_1/\log n\leq \inf_{\theta_0\in\Theta^{\beta}(M)}\underline\a_n(\theta_0)
\leq\sup_{\theta_0\in\Theta_s^{\beta}(M)}\overline\a_n(\theta_0)
\leq \beta+K_2/\log n,\label{eq: BoundsforBounds}
\end{align}
for some positive constants $K_1$ and $K_2$ depending only on $\beta_{\min},\beta_{\max},M,C,\rho$ and $\eps$.

Then for convenience we introduce the notation
\begin{align}
r_{n,\gamma}^2(\a)&=D^2\zeta_{\gamma/2}^2s_n^2(\a)=D^2\zeta_{\gamma/2}^2\sum_{i=1}^{\infty}\frac{l_i^2\kappa_i^{-2}}{i^{1+2\a}\kappa_i^{-2}+n}.\label{def: r}
\end{align}

Using the notations of Section \ref{Proof: counter}
the coverage of the empirical Bayes credible set, similarly to the inequality $\eqref{eq: coverCounter}$, can be bounded from below by
\begin{align}
&\inf_{\theta_0\in\Theta_s(M)}\P_{\theta_0}\Big(L\th_0\in \hat{C}_n^{L}(D),\,\forall L\in  L_s^{q}(R) \Big)\label{eq: BiasLin}\\
&\geq
\inf_{\theta_0\in\Theta_s(M)}\PP_{\theta_0}\Big(B_n^{L}(\hat\alpha_n-1/2)+V_n(\hat\alpha_n-1/2)\leq  r_{n,\gamma}(\hat\alpha_n-1/2),\,\forall L\in L_s^{q}(R)\Big)\nonumber\\
&\geq \inf_{\theta_0\in\Theta_s(M)}\P_{\theta_0}\Big(\sup_{\substack{\alpha\in [\underline\alpha_n-1/2,\overline\alpha_n-1/2] \\
L\in L_s^{q}(R)}}V_n(\a)  \leq \inf_{\substack{\alpha\in [\underline\alpha_n-1/2,\overline\alpha_n-1/2] \\
L\in L_s^{q}(R)}}\big[r_{n,\gamma}(\alpha) -B_n^{L}(\alpha)\big]\Big)-o(1).\nonumber
\end{align}

Therefore it is sufficient to show that there exist constants $C_1,C_2$ and $C_3$ satisfying $C_3>C_2+C_1$ such that for all $\beta\in[\beta_{\min},\beta_{\max}]$ and $\theta_0\in \Theta_s^{\beta}(M)$

\begin{align}
\sup_{\a\in [\underline\alpha_n-1/2,\overline\alpha_n-1/2]}B_n^{L}(\alpha)\leq C_1 n^{-\frac{\beta+q}{2\beta+2p}},\label{eq: 3eq1}\\
P_{\theta_0}\Big(\sup_{\a\in [\underline\alpha_n-1/2,\overline\alpha_n-1/2]}V_n(\alpha)\leq C_2 n^{-\frac{\beta+q}{2\beta+2p}}\Big)\rightarrow1,\label{eq: 3eq2}\\
\inf_{\a\in [\underline\alpha_n-1/2,\overline\alpha_n-1/2]}r_{n,\gamma}(\alpha)\geq C_3 n^{-\frac{\beta+q}{2\beta+2p}}.\label{eq: 3eq3}
\end{align}

We first deal with inequality $\eqref{eq: 3eq3}$. Applying assumptions $\eqref{assump: LinFuncSelf}$ and $\eqref{def: r}$ one can obtain that
\begin{align*}
r_{n,\gamma}^2(\a)&\geq \frac{D^2\zeta_{\gamma/2}^2}{C^{2}}\sum_{i=j_0/K+1}^{\infty}\frac{(K i)^{2p}}{((K+1)i)^{1+2\a+2p}+n}\sum_{j=iK}^{(i+1)K-1}l_j^2\nonumber\\
&\geq
\frac{D^2\zeta_{\gamma/2}^2K^{2p-2q}}{R^2C^{2}(K+1)^{1+2\a+2p}}\sum_{i=j_0/K+1}^{\infty}\frac{i^{-1-2q+2p}}{i^{1+2\a+2p}+n},
\end{align*}
which following from Lemma \ref{Lemma: TechLem10} is further bounded from below by constant times $D^2n^{-\frac{1+2\a+2q}{1+2\a+2p}}$ for $1+2\a+2q>0$ and infinity else. Therefore by applying the inequality $\eqref{eq: BoundsforBounds}$ we get that $C_3$ can be arbitrary large for a large enough choice of $D^2$.

Next we deal with the convergence $\eqref{eq: 3eq2}$. From the proof of Theorem 5.3 of \cite{Knapik} we get that $V_{n}(\a)=|t_n(\a)Z|$ with $Z$ a standard normal random variable and
\begin{align*}
t_n^2(\a)= \sum_i^{\infty}\frac{nl_i^2 \kappa_i^{-2}}{(i^{1+2\a}\kappa_i^{-2}+n)^2}\leq R^2C^{4}\sum_i^{\infty}\frac{n i^{-1-2q+2p}}{(i^{1+2\a+2p}+n)^2}.
\end{align*}
The right hand side of the preceding display similarly to $s_n^2(\a)$ is bounded above by constant times $n^{-\frac{1+2\a+2q}{1+2\a+2p}}$ for $1+2\a+2q>0$ and infinity otherwise. Then following from the inequality $\eqref{eq: BoundsforBounds}$ one can obtain for $q<p$ that
\begin{align*}
\sup_{\a\in[\underline\a_n-1/2,\overline\a_n-1/2]}t_n(\a)&=t_n(\underline\a_n-1/2)
\lesssim n^{-\frac{\underline\a_n+q}{2\underline\a_n+2p}}
\lesssim n^{-\frac{\beta+q}{2\beta+2p}},
\end{align*}
providing us the convergence$\eqref{eq: 3eq2}$.

Finally we deal with the bias term $\eqref{eq: 3eq1}$. Following from assumptions $\eqref{def: hyperrectangle}$ and $\eqref{assump: LinFuncSelf}$ we have
\begin{align*}
|B_n^{L}(\a)|\leq\sum_{i=1}^{\infty}\frac{|l_i \theta_{0,i}|i^{1+2\a}\kappa_i^{-2}}{i^{1+2\a}\kappa_i^{-2}+n}
\leq C^2RK \sum_{i=1}^{\infty}\frac{i^{2\a+2p-\beta-q}}{i^{1+2\a+2p}+n}.
\end{align*}
From the inequality $\eqref{eq: BoundsforBounds}$ we have for $\a\geq\underline\a_n-1/2$ and large enough $n$ that the inequality $\beta+q<1+2\a+2p$ holds, hence the preceding inequality is further bounded from above by constant times $n^{-\frac{\beta+q}{1+2\a+2p}}$ by applying Lemma \ref{Lemma: TechLem10} (with $m=0$). So we can conclude that for $\a\geq\underline\a_n-1/2\geq \beta-1/2-K_1/\log n$
\begin{align*}
|B_n^L(\a)|\lesssim n^{-(\beta+q)/(1+2\a+2p)}\lesssim n^{-(\beta+q)/(2\beta+2p)}.
\end{align*}

To prove adaptivity we note that following again from the inequality $\eqref{eq: BoundsforBounds}$ we have
\begin{align*}
\sup_{\a\in[\underline\a_n-1/2,\overline\a_n-1/2]} s_n(\a)
\lesssim n^{-(\underline\a_n+q)(2\underline\a_n+2p)}\lesssim n^{-(\beta+q)/(2\beta+2p)}.
\end{align*}

\vspace{-0.5cm}
\section{Appendix}
\begin{lemma}[Lemma 10.2 of \cite{SzVZ2}]\label{Lemma: TechLem10}
For any $l, m, r, s\ge 0$ with $c:=lr-s-1>0$ and $n \ge e^{(2mr/c)\vee r}$,
$$(3^r+1)^{-l}\bigl({\log n}/r\bigr)^{m}n^{-c/r}\le \sum_{i=1}^\infty \frac{i^{s}(\log i)^m}{(i^{r}+n)^l}
\le  (3+2c^{-1})\bigl({\log n}/r\bigr)^{m}n^{-c/r}.$$
\end{lemma}
%\begin{acknowledgement}
%The author would like to thank Bartek Knapik for reading the manuscript and suggesting changes, which improved the quality of the manuscript.
%\end{acknowledgement}
\begin{acknowledgement}
The author would like to thank Harry van Zanten, Bartek Knapik and an anonymous referee for their suggestions leading to better readability and presentation of the manuscript. 
\end{acknowledgement}
\vspace{-0.7cm}
\bibliographystyle{acm}
\bibliography{references}
\end{document}